\documentstyle[amssymb,12pt]{article}

\begin{document}

\newtheorem{theorem}{Theorem}{}
\newtheorem{lemma}[theorem]{Lemma}{}
\newtheorem{corollary}[theorem]{Corollary}{}
\newtheorem{conjecture}[theorem]{Conjecture}{}
\newtheorem{proposition}[theorem]{Proposition}{}
\newtheorem{axiom}{Axiom}{}
\newtheorem{remark}{Remark}{}
\newtheorem{example}{Example}{}
\newtheorem{exercise}{Exercise}{}
\newtheorem{definition}{Definition}{}

\title{Reduction of quantum systems with arbitrary first class constraints and 
Hecke algebras}

\author{
Alexey Sevostyanov \footnote{e-mail seva@teorfys.uu.se}\\ %
Institute of Theoretical Physics, Uppsala University,\\ and Steklov Mathematical 
Institute, St.Petersburg 
}

\maketitle
\begin{flushright}
UU--ITP 5/98

\end{flushright}

\begin{abstract}
We propose a method for reduction of quantum systems with arbitrary first class 
constraints. 
An appropriate mathematical setting for the problem is homology of associative 
algebras.
For every such an algebra $A$ and its subalgebra $B$ with an augmentation 
$\varepsilon$
there exists a cohomological complex which is 
a generalization of the BRST one. Its cohomology is an associative graded 
algebra $Hk^{*}(A,B)$ 
which we call 
the Hecke algebra of the triple $(A,B,\varepsilon)$. It acts in the cohomology 
space $H^{*}(B,V)$ for 
every left $A$--
module $V$. In particular the zeroth graded component $Hk^{0}(A,B)$ acts in the 
space of $B$--
invariants of $V$ and provides the reduction of the quantum system.  
\end{abstract}

\section*{Introduction}

 The purpose of this paper is to generalize the well--known BRST quantization 
procedure to arbitrary
associative algebras. An appropriate mathematical setting for the BRST 
cohomology of Lie algebras 
\cite{BRS} was
proposed by Kostant and Sternberg in \cite{KSt}. Their observation was that 
given a Hamiltonian action of 
a Lie group on a Poisson manifold one can construct a super--Poisson complex 
whose zeroth cohomology is
the algebra of functions on the reduced space over zero value of the 
corresponding moment map \cite{A} , \cite{MW}. 
This complex 
admits a quantization and the zeroth cohomology of the quantum complex may be 
treated as a quantization of 
the classical reduced space. In that case the quantum counterparts of matrix 
elements of the moment map 
form a Lie algebra and represent a system of the first--class constraints for 
the quantum reduction.

Our approach to the BRST cohomology differs from the one described above. We 
start with the quantum complex 
directly. It turns out that the quantum BRST cohomology may be defined using the 
language of homological 
algebra , resolutions , etc.. The particular complex proposed in \cite{KSt} 
corresponds to the 
standard resolution of the ground field. As usual , in homological algebra 
different choices of resolutions 
lead to the same homology. This allows us to generalize the BRST cohomology to 
arbitrary systems of 
the first--class constraints.

The BRST quantization is not unique in the following sense. One can always 
define the BRST cohomology 
related to the usual cohomology of representation spaces of the quantum system. 
However under some 
technique assumptions there exists another version of the BRST reduction related 
to the semiinfinite 
cohomology of the representation spaces \cite{F}. In the case of Lie algebras it 
requires a normal ordering in the 
differential \cite{KSt}. In the present paper we only discuss the usual BRST 
cohomology. The latter one will be 
explained in a subsequent paper.

\section{Endomorphisms of complexes}\label{end}

Let $A$ be an associative ring with unit , $X$ be a graded complex of left $A$ 
-- modules equipped 
with a differential of degree -1. Recall the definition of the complex 
$Y=End_{A}(X)$ \cite{MacLane}.

By the definition $Y$ is a $\Bbb Z$ graded complex

\begin{equation}
Y=\bigoplus_{n=-\infty}^{\infty} Y^n
\end{equation}

with graded components defined as

\begin{equation}
Y^n=\prod_{p+q=n} Y^{p,q},
\end{equation}

where 

\begin{equation}
Y^{p,q}=Hom_{A}(X^{p},X^{-q}).
\end{equation}

Clearly $Y$ is a subalgebra in the full algebra of $A$--endomorphisms of $X$. It 
is easy to see that $Y$ is
closed with respect to the multiplication given by composition of endomorphisms. 
Thus it is a graded 
associative algebra.
We emphasize that $Y$ is not a bigraded space. 

Introduce a differential on $Y$ of degree +1 as follows

\begin{equation}
\begin{array}{c}
({\bf d}f)^{p,q}=(-1)^{p+q}f{p-1,q}\circ d +d\circ f^{p,q-1}, \\
f=\{ f^{p,q} \}, f^{p,q} \in Y^{p,q}.
\end{array}
\end{equation}

If $f$ is homogeneous then

\begin{equation}\label{diff}
{\bf d}f=d\circ f -(-1)^{deg (f)}f\circ d .
\end{equation}

So that ${\bf d}$ is the supercommutator by $d$

We shall consider also the partial differentials $d'$ and $d''$ on Y :

\begin{equation}\label{part}
\begin{array}{cc}
\mbox{for } f\in Y^{p,q}& \\
(d'f)(x)=(-1)^{p+q+1}f(dx),& x\in X^{p+1}; \\
(d''f)(x)=df(x),& x\in X^{p}.
\end{array}
\end{equation}

It is easy to check that

\begin{equation}
d'^2=d''^2=d'd''+d''d'=0
\end{equation}

These conditions ensure that ${\bf d}^2=0$.

The following property of ${\bf d}$ is crutial for the subsequent 
considerations.

\begin{lemma}
${\bf d}$ is a superderivation of $Y$.
\end{lemma}

{\em Proof.} Let $f$ and $g$ be homogeneous elements of $Y$. Then $deg (fg) = 
deg (f) + deg (g)$ and 
(\ref{diff}) yields:

\begin{equation}
\begin{array}{l}
{\bf d}(fg)=d\circ fg - (-1)^{deg (f) +deg (g)} fg \circ d= \\ 
d\circ fg -(-1)^{deg (f)} f \circ d \circ g +(-1)^{deg (f)} f \circ d \circ g - 
(-1)^{deg (f) +deg (g)} fg \circ d= \\
({\bf d}f)g+(-1)^{deg (f)}f({\bf d}g). 
\end{array}
\end{equation}

The proof follows.

The most important consequence of the lemma is

\begin{theorem}\label{alg}
The homology space $H^*(Y)$ inherits a multiplicative structure from Y. 
Thus $H^*(Y)$ is a graded associative algebra. 
\end{theorem}

{\em Proof.} First, the product of two cocycles is a cocycle. 
For if $f$ and $g$ are homogeneous and ${\bf d}f={\bf d}g=0$ then

\begin{equation}
{\bf d}fg=({\bf d}f)g+(-1)^{deg (f)}f({\bf d}g)=0.
\end{equation}

Now we have to show that the product of homology classes is well defined. 
It is suffices to verify that the product of a homogeneous cocycle with a 
homogeneous cobondary is homologous to zero.
For instance consider the product $f{\bf d}h$. Then (\ref{diff}) gives

\begin{eqnarray}
f{\bf d}h=f\circ (d\circ h -(-1)^{deg (h)}h\circ d)= \\
(-1)^{deg (f)}d\circ f \circ h - (-1)^{deg (h)}f\circ h\circ d = 
(-1)^{deg(f)}{\bf d}(fh). \nonumber
\end{eqnarray}

This completes the proof.

One of the principal statements of homological algebra says that homotopically 
equivalent complexes have
the same homology. In particular the vector space $H^*(Y)$ depends only on the  
homotopy
class of the complex $X$. It turns out that the same is true for the algebraic 
structure of $H^*(Y)$.
Indeed we have the following

\begin{theorem}\label{equiv}
Let $X, X'$ be two homotopically equivalent graded complexes of left 
$A$--modules. Then

\begin{equation}
 H^*(Y)\simeq H^*(Y')
\end{equation}

as graded associative algebras.
\end{theorem}

{\em Proof.} Let $F:X \rightarrow X' , F':X' \rightarrow X$ be two maps of the 
complexes such that

\begin{equation}
\begin{array}{ll}
F'F-id_X=d_Xs+sd_X,& s:X\rightarrow X ,\\
FF'-id_{X'}=d_{X'}s'+s'd_{X'},& s:X'\rightarrow X' ,\\
s \in Y^{-1},& s' \in Y'^{-1}. 
\end{array}
\end{equation}

Consider the induced mappings of the complexes $Y,Y'$:

\begin{equation}
\begin{array}{c}
FF'^{*}:Y \rightarrow Y' ,\\
FF'^{*}f= F \circ f \circ F' , f\in Y ;\\
F'F^{*}:Y' \rightarrow Y ,\\
F'F^{*}g= F' \circ g \circ F , f\in Y' .\\
\end{array}
\end{equation}

Their compositions are homotopic to the identity maps of $Y$ and $Y'$ (see 
\cite{carteil} , Chap.4 for a 
general statement about equivalences of functors). But it means that $FF'^{*}$ 
is inverse to $F'F^{*}$ 
when restricted to homology. Thus $H^*(Y)$ is isomorphic to $H^*(Y')$ as a 
vector space. We have to
show that the restrictions of $FF'^{*}$ and $F'F^{*}$ to the homologies are 
homomorphisms of algebras.

Let $f$ and $g$ be homogeneous elements of $Y$ and ${\bf d}_Xf={\bf d}_Xg=0$. By 
the definition of the
induced maps we have

\begin{equation}
FF'^{*}(fg)=F \circ fg \circ F'.
\end{equation}

On the other hand

\begin{eqnarray}\label{hom}
FF'^{*}(f)FF'^{*}(g)=F \circ f \circ F'F \circ g \circ F'=   \\ 
F \circ f(id_X+d_Xs+sd_X)g \circ F'.\nonumber
\end{eqnarray}

Now recall that $f$ and $g$ are cocycles in $Y$. By (\ref{diff}) they 
supercommute with $d_X$:

\begin{equation}\label{cocycle}
d_X\circ f =(-1)^{deg (f)}f\circ d_X.
\end{equation}

Using (\ref{cocycle}) and the fact that $F$ and $F'$ are morphisms of complexes 
we can rewrite (\ref{hom})
as follows

\begin{eqnarray}\label{hom1}
F \circ f(id_X+d_Xs+sd_X)g \circ F'= \nonumber \\
F \circ fg \circ F' + (-1)^{deg f}d_{X'} \circ F \circ fsg \circ F' + 
(-1)^{deg g}F \circ fsg \circ F' \circ d_{X'} = \\
F \circ fg \circ F'+(-1)^{deg f}{\bf d}_{X'}(F \circ fsg \circ F') . \nonumber
\end{eqnarray}

Finally observe that by (\ref{hom1}) $FF'^{*}(fg)$ and $FF'^{*}(f)FF'^{*}(g)$ 
belong to the same 
homology class in $H^{*}(Y')$. This completes the proof.

\section{Hecke algebras}\label{Hecke}

Let $A$ be an associative algebra over a ring $K$ with unit , $B$ be its 
subalgebra with an
augmentation $\varepsilon : B\rightarrow K$ ($\varepsilon $  is a homomorphism 
of $K$--algebras).

Let $X$ be a projective resolution of the left $B$--module $K$. Then the complex

\begin{equation}
A \otimes_B X
\end{equation}

has the natural structure of a left $A$--module. We can apply theorem \ref{alg} 
to define a graded associative 
algebra

\begin{equation}
Hk^*(A,B)=H^*(End_A(A \otimes_B X)).
\end{equation}

Observe that all $B$--projective resolutions of $K$ are homotopically 
equivalent. Hence
by theorem \ref{equiv} $Hk^*(A,B)$ does not depend on the resolution $X$. We 
shall call it the Hecke 
algebra of the the triple $(A,B,\varepsilon )$.

Now consider $A$ as a left $A$--module and a right $B$--module via 
multiplication. In this way $A$ becomes
an $A\otimes B^{opp}$--left module. Let $X'$ be a projective resolution of the 
module. The complex

\begin{equation}
X' \otimes_B K
\end{equation}

is a left $A$--module. Therefore there exists an associative algebra

\begin{equation}
\widehat{Hk}^*(A,B)=H^*(End_A(X' \otimes_B K))
\end{equation}

independent of the resolution $X'$.

\begin{theorem}\label{iso}
$Hk^*(A,B)$ is isomorphic to $\widehat{Hk}^*(A,B)$ as a graded associative 
algebra.
\end{theorem}

{\em Proof.} We shall use the standard bar resolutions for computation of 
$\widehat{Hk}^*(A,B)$ and 
${Hk}^*(A,B)$ \cite{MacLane} , \cite{carteil}. 
Consider the complex $B\otimes T(I(B)) \otimes B$, where $I(B)=B/K$ and $T$ 
denotes the tensor algebra 
of the vector space. Elements of $B\otimes T(I(B)) \otimes B$ are usually 
written as
$a[a_1,\ldots ,a_s]a'$. The differential is given by

\begin{eqnarray}
da[a_1,\ldots ,a_s]a'=aa_1[a_2,\ldots ,a_s]a'+ \\
\sum_{k=1}^{s-1}(-1)^{k}a[a_1,\ldots ,a_ka_{k+1},\ldots  ,a_s]a' + 
(-1)^sa[a_1,\ldots ,a_{s-1}]a_sa'. \nonumber
\end{eqnarray}

Then $B\otimes T(I(B)) \otimes B \otimes_B K = B\otimes T(I(B)) \otimes K$ is a 
free resolution of
the left $B$--module $K$. And $A \otimes_A B\otimes T(I(B)) \otimes B =A\otimes 
T(I(B)) \otimes B$ is 
a free resolution of $A$ as a right $B$--module. The complex $A\otimes T(I(B)) 
\otimes B$ is also a 
left free $A$--module via the left multiplication by elements from $A$. Hence 
this is an $A\otimes B^{opp}$--
free resolution of $A$.  

Thus the complex $End_A(A\otimes_B B\otimes T(I(B)) \otimes K)=End_A(A\otimes 
T(I(B)) \otimes K)$ 
for computation of $Hk^*(A,B)$ is canonically isomorphic to the complex
$End_A(A \otimes T(I(B)) \otimes B \otimes_B K)=End_A(A \otimes T(I(B)) \otimes 
K)$ for computation
of $\widehat{Hk}^*(A,B)$. This establishes the isomorphism of the algebras.

\section{Action in homology and cohomology spaces}

Recall that for every left $B$--module $V$ the cohomology modules are defined to 
be

\begin{equation}\label{cohomol}
H^*(V)=Ext_B^*(K,V)=H^*(Hom_B(X,V)),
\end{equation}

where $X$ is a projective resolution of $K$. While for every right $B$--module 
$W$ one can define
the homology modules

\begin{equation}\label{homol}
H_*(W)=Tor_*^B(W,K)=H_*(W\otimes_B X).
\end{equation}

Now observe that for every right $A$--module $V$ 
the complex (\ref{cohomol}) for calculation its cohomology as a right 
$B$--module 
may be represented as follows:

\begin{equation}\label{cohomcompl}
Hom_B(X,V)=Hom_A(A\otimes_B X ,V).
\end{equation}

Endow the space $Hom_A(A\otimes_B X ,V)$ with a right $End_A(A\otimes_B 
X)$--action:

\begin{equation}\label{act1}
\begin{array}{ll}
Hom_A(A\otimes_B X ,V) \times End_A(A\otimes_B X) \rightarrow Hom_A(A\otimes_B X 
,V) , &\\
\varphi \times f \mapsto  \varphi \circ f, &\\
\varphi \in Hom_A(A\otimes_B X ,V), f \in End_A(A\otimes_B X).&
\end{array}
\end{equation}

The action is well defined since $f$ commutes with the left $A$--action. Clearly 
this action respects the
gradings, i.e., it is an action of the graded associative algebra on the graded 
module.

\begin{theorem}\label{cohomact}

The action (\ref{act1}) gives rise to a right action

\begin{eqnarray}\label{cohomact1}
H^*(V)\times {Hk}^*(A,B) \rightarrow H^*(V), \\
H^n(V)\times {Hk}^m(A,B) \rightarrow H^{n+m}(V).\nonumber
\end{eqnarray}

\end{theorem}

{\em Proof.} Let $\varphi \in Hom_A(A\otimes_B X ,V)$ and $d\varphi =\varphi 
\circ d =0$. Let also 
$f \in End_A(A\otimes_B X)$ be a homogeneous cocycle. By (\ref{cocycle}) 
$\varphi \circ f$ is a cocycle
in $Hom_A(A\otimes_B X ,V)$. Indeed

\begin{equation}
d(\varphi \circ f)=\varphi \circ f \circ d =(-1)^{deg (f)} \varphi \circ d \circ 
f =0.
\end{equation}

Then we need to show that the action does not depend on the choice of the 
representative $f$ in the
homology class $[f]$,that is $\varphi \circ {\bf d}g$ is homologous to zero for 
every homogeneous 
$g\in End_A(A\otimes_B X)$. This is a direct consequence of the definitions:

\begin{equation}
\varphi \circ {\bf d}g= \varphi \circ (d\circ g - (-1)^{deg(g)} g\circ d)= 
-(-1)^{deg(g)}d(\varphi \circ g),
\end{equation}

since $\varphi \circ d =0$. Finally let us check that the action is independent 
of the representative
in the homology class $[\varphi]$. For $\psi \in Hom_A(A\otimes_B X ,V)$ $d\psi 
\circ f$ is always 
homologous to zero:

\begin{equation}
d\psi \circ f= \psi \circ d \circ f = (-1)^{deg(f)}\psi \circ f \circ d 
=(-1)^{deg(f)}d(\psi \circ f).
\end{equation}

This concludes the proof.

Similarly for every right $A$--module $W$ one can equip the homology module 
$H_*(W)$ with a 
structure of left ${Hk}^*(A,B)$--module. First the complex $W\otimes_B 
X=W\otimes_A A\otimes_B X$ has
a natural structure of left $End_A(A\otimes_B X)$--module:

\begin{equation}\label{act2}
\begin{array}{ll}
End_A(A\otimes_B X)\times W\otimes_A A\otimes_B X  \rightarrow W\otimes_A 
A\otimes_B X , &\\
f \times w\otimes x\mapsto  w\otimes f(x), &\\
w\otimes x \in W\otimes_A (A\otimes_B X), f \in End_A(A\otimes_B X).&
\end{array}
\end{equation}

Observe that according to the convention of section \ref{end} elements of 
$End_A^n(A\otimes_B X)$ have
degree -1 as operators in the graded space $W\otimes_A A\otimes_B X$:

\begin{equation}
End_A^n(A\otimes_B X)\times W\otimes_A A\otimes_B X_m \rightarrow W\otimes_A 
A\otimes_B X_{m-n}.
\end{equation}

The following assertion is an analogue of theorem \ref{cohomact} for homology.

\begin{theorem}
The action (\ref{act2}) gives rise to a left action

\begin{eqnarray} 
Hk(A,B)^* \times H_*(W) \rightarrow H_*(W),\\
Hk(A,B)^n \times H_m(W) \rightarrow H_{m-n}(W).\nonumber
\end{eqnarray}

\end{theorem}

\section{Structure of the Hecke algebras}

In this section we investigate the Hecke algebras under some technique 
assumptions. The main theorem here is

\begin{theorem}\label{struct}
Assume that

\begin{equation}
Tor_n^B(A,K)=0 \mbox{  for } n>0.
\end{equation}

Then

\begin{equation}
Hk^n(A,B)=Ext^n_A(A\otimes_BK,A\otimes_BK)=Ext^n_B(K,A\otimes_BK).
\end{equation}

In particular

\begin{equation}
\begin{array}{l}
Hk^n(A,B)=0 , n<0;\\
Hk^0(A,B)=Hom_B(K,A\otimes_BK).
\end{array}
\end{equation}

\end{theorem}

{\em Proof.} Equip the complex $Y=End_A(A\otimes T(I(B)) \otimes K)$ , which we 
used in theorem \ref{iso} for
computation of $Hk^*(A,B)$ , with the first filtration as follows:

\begin{equation}
F^kY=\sum_{n=-\infty}^{\infty}\prod_{p+q=n , p\geq k}Y^{p,q}.
\end{equation}

The associated graded complex with respect to the filtration is the double 
direct sum

\begin{equation}
GrY=\sum_{p,q=-\infty}^{\infty}Y^{p,q}.
\end{equation}

One can show that the filtration is regular and the second term of the 
corresponding spectral sequence is

\begin{equation}\label{spec}
E_2^{p,q}=H^p_{d'}(H^q_{d''}(GrY)),
\end{equation}

where $H^*_{d'}$ and $H^*_{d''}$ denote the homologies of the complex with 
respect to the partial differentials (\ref{part}).

 Now observe that at the same time the complex $A\otimes T(I(B)) \otimes K$ is a 
complex 
for calculation
of $Tor_n^B(A,K)$ because $A\otimes T(I(B)) \otimes B$ is a free resolution of 
$A$ as a right $B$--module. It is
also free as a left $A$--module. Therefore the functor $Hom_A(A\otimes T(I(B)) 
\otimes K, \cdot )$ is exact.
By the assumption $H^*(A\otimes T(I(B)) \otimes K)=Tor_0^B(A,K)=A\otimes_BK$. 
Using the last two
observations we can calculate the cohomology of the complex $GrY$ with respect 
to the differential $d''$ :

\begin{equation}\label{degener}
\begin{array}{l}
H^*_{d''}(GrY)=H^*_{d''}(Hom_A(A\otimes T(I(B)) \otimes K,A\otimes T(I(B)) 
\otimes K)=\\
Hom_A(A\otimes T(I(B)) \otimes K,A\otimes_BK).
\end{array}
\end{equation}

Here $Hom_A$ should be thought of as the direct sum of the double graded 
components. Now (\ref{degener}) provides that 
the spectral sequence (\ref{spec}) degenerates at the second term. Moreover

\begin{equation}
E_2^{p,*}=H^p_{d'}(H^0_{d''}(GrY))=H^p_{d'}(A\otimes T(I(B)) \otimes 
K,A\otimes_BK).
\end{equation}

But the complex $A\otimes T(I(B)) \otimes K$ may be regarded as a free 
resolution of the left 
$A$--module $A\otimes_BK$. Therefore

\begin{equation}
E_2^{p,*}=Ext^p_A(A\otimes_BK,A\otimes_BK).
\end{equation}

Finally by theorem 5.12 \cite{carteil} we have:

\begin{equation}
Hk^n(A,B)=H^n(Y)=E_2^{n,0}=Ext^n_A(A\otimes_BK,A\otimes_BK).
\end{equation}

Since $Tor_n^B(A,K)=0 \mbox{  for } n>0$ we can apply Shapiro lemma to simplify 
the last expression:

\begin{equation}
Ext^n_A(A\otimes_BK,A\otimes_BK)=Ext^n_B(K,A\otimes_BK).
\end{equation}

This completes the proof.

\begin{remark}
In particular the conditions of the theorem are satisfied if $A$ is projective 
as a right $B$--module. For instance
suppose that there exists a subspace $N \subset A$ such that the multiplication 
in $A$ provides an isomorphism of 
the vector spaces $A \simeq N\otimes B$. Then $A$ is a free right $B$--module.
\end{remark}

\section{Comparison with the BRST complex}

Let $\frak g$ be a Lie algebra over a field $K$. For simplicity we suppose that 
$\frak g$ is finite--dimensional. However 
the aguements presented below remain to be true, with some technique 
modifications, for an arbitrary Lie algebra.
We shall apply the construction of section \ref{Hecke} in the following 
situation. 

Let $A$ be an associative 
algebra over $K$ and $B=U({\frak g})$ be its subalgebra. Note that $U({\frak 
g})$ is naturally augmented.
Consider the $U({\frak g})$--free resolution of the left $U({\frak g})$--module 
$K$ as follows:

\begin{equation}
\begin{array}{l}
X=U({\frak g})\otimes \Lambda ({\frak g}),\\
d (u\otimes x_1 \wedge \ldots \wedge x_n)= 
\sum_{i=1}^n (-1)^{i+1} ux_i\otimes x_1 \wedge \ldots \wedge \widehat{x_i} 
\wedge \ldots \wedge x_n +\\
\sum_{1\leq i< j \leq n}(-1)^{i+j} u\otimes [x_i,x_j] \wedge 
x_1 \wedge \ldots \wedge \widehat{x_i} \wedge \ldots \wedge \widehat{x_j} \wedge 
\ldots \wedge x_n,
\end{array}
\end{equation}

where the symbol $\widehat{x_i}$ indicates that $x_i$ is to be omitted. Then

\begin{equation}
A\otimes_{U({\frak g})} X \simeq A\otimes \Lambda ({\frak g})
\end{equation}

is a complex with a differential given by the operator

\begin{equation}\label{different}
d=\sum_i e_i\otimes e_i^*-\sum_{i,j} 1\otimes [e_i,e_j]e_i^*e_j^*,
\end{equation}

here $e_i$ is a linear basis of ${\frak g}$, $e_i^*$ is the dual basis; 
$e_i\otimes 1$ being regarded as the 
operator of right multiplication in $A$ and $1\otimes e_i , 1\otimes e_i^* $ as 
the operators of exterior and 
inner multiplications in $\Lambda ({\frak g})$ respectively. 

Now observe that

\begin{equation}
\begin{array}{l}
End_A(A\otimes \Lambda ({\frak g}))=Hom_K(\Lambda ({\frak g}),A^{opp}\otimes 
\Lambda ({\frak g}))=\\
A^{opp}\otimes End_K(\Lambda ({\frak g}))=A^{opp} \otimes C({\frak g}+{\frak 
g}^*),
\end{array}
\end{equation}

where $C({\frak g}+{\frak g}^*)$ is the Clifford algebra of the space ${\frak 
g}+{\frak g}^*$. Under this 
identification $A^{opp}$ acts on $A\otimes \Lambda ({\frak g})$ by the 
multiplications in $A$ from the right and the Clifford 
algebra acts by the exterior and inner multiplications in $\Lambda ({\frak g})$. 
This allows to consider the 
differential (\ref{different}) as an element of the complex $A^{opp} \otimes 
C({\frak g}+{\frak g}^*)$. 

It is easy to see that the canonical $\Bbb Z$ grading of the complex $A^{opp} 
\otimes C({\frak g}+{\frak g}^*)$ coincides 
mod 2 with the ${\Bbb Z}_2$ grading inherited from the Clifford algebra. 
Therefore according to (\ref{diff}) 
the differential $\bf d$ is given by the supercommutator in $A^{opp} \otimes 
C({\frak g}+{\frak g}^*)$ by element 
(\ref{different}). This establishes

\begin{theorem}
The complex $(End_A(A\otimes_{U({\frak g})} X) , {\bf d})$ is isomorphic to the 
BRST one 
$A^{opp} \otimes C({\frak g}+{\frak g}^*)$ with the differential being the 
supercommutator by element (\ref{different}).
\end{theorem}

\section{Relation to the quantum reduction}

The results of the previous section imply that if $K$ is the field of complex 
numbers $\Bbb C$ then 
$Hk^0(A,U({\frak g}))^{opp}$ may be thought of as a result of the quantum 
reduction in $A$ with $U({\frak g})$ being a system of the first--class 
constraints \cite{KSt}. We shall show that 
this treatment remains to be true in the general situation of section 
\ref{Hecke}.

Suppose that $A$ is a quantization of a classical system, that is , $A$ is 
included into a family of 
associative algebras $A_h$ parametrized by a complex number $h$ such that for 
different $h$ 
$A_h$ are isomorphic as vector spaces , $A_0$ is commutative and the formula

\[
\{ a,b \} =\lim_{h\rightarrow 0} {ab - ba \over h}
\]

defines a Poisson algebra structure on $A_0$. The classical limit of $B$ is a 
Poisson subalgebra $B_0$ in $A_0$ 
with a character $\varepsilon _0 : B_0 \rightarrow {\Bbb C}$. Let $J_0$ be the 
ideal in $A_0$ generated by the 
kernel of the map $\varepsilon _0$. Then the classical reduced Poisson algebra 
coincides with the subspace 
of Poisson $B_0$--invariants in the quotient $A_0 / J_0$. In typical situations 
$A_0$ is the Poisson algebra 
of functions on a Poisson manifold. In this case the scheme of the reduction was 
suggested by Direc in 
\cite{D}.

Now assume that the conditions of theorem \ref{struct} are satisfied. Then the 
algebra $Hk^0(A,B)^{opp}$ is 
isomorphic to the algebra of $B$--invariants in the quotient $A/J$ where $J$ is 
the left ideal in $A$ 
generated by the kernel of the augmentation map $\varepsilon$. Thus the 
classical limit of $Hk^0(A,B)^{opp}$ 
is exactly the reduced Poisson algebra defined above.

If theorem \ref{struct} does not hold the algebra $Hk^0(A,B)^{opp}$ can be still 
treated as a quantization 
of the classical reduced space in the following sense.

Recall the scheme of the Dirac quantum reduction \cite{D}. Suppose again that we 
are given a quantum system with first--class 
constraints, that is an associative algebra $A$ over $\Bbb C$ together with a 
representation $V$ and its subalgebra 
$B$ equipped with a character $\varepsilon : B \rightarrow {\Bbb C}$. According 
to Dirac the space of the physical states 
for the reduced system is the space of $B$--invariants in $V$:

\begin{equation}
V^B=\{ v \in V : bv=\varepsilon (b)v \mbox{ for every } b\in B \},
\end{equation}

that is the zeroth cohomology space of $V$ as a left $B$--module 
$H^0(V)=V^B=Hom_B({\Bbb C},V)$.
And the algebra of observables $A^B_V$ of the reduced system is formed by 
operators $a \in A$ such that

\begin{equation}\label{dirac}
[b,a]v=0 \mbox{ for every } b \in B , v \in V^B.
\end{equation}

This condition is equivalent to

\begin{equation}
bav=\varepsilon (b)av.
\end{equation}

It means that the space $V^B$ is invariant with respect to the action of 
$A_V^N$. Clearly, such operators form an 
algebra. 

From the other side we have an action (\ref{cohomact1}) of the algebra 
$Hk^0(A,B)^{opp}$ in the space $V^B$. 
From the definition of the action it is clear that 
only $(0,0)$--bidegree components give a nontrivial contribution to the action. 
They may be represented by 
elements from $A$. Moreover

\begin{theorem}
The algebra $Hk^0(A,B)^{opp}$ satisfies condition (\ref{dirac}) for every 
representation $V$ of $A$. So it may
be regarded as a universal Dirac reduction of the physical system.
\end{theorem}

The statement of the theorem follows directly from theorem \ref{cohomact}. 
Namely , condition (\ref{dirac}) 
ensures that the action of the algebra $Hk^0(A,B)^{opp}$ in the space of 
invariants is well--defined.

\end{document}